\pdfoutput=1
\documentclass[10pt]{article}
\usepackage[a4paper, margin=1in]{geometry}
\usepackage{amsmath, amssymb, amsthm}
\usepackage[utf8]{inputenc}
\usepackage[english]{babel}
\usepackage[numbers]{natbib}
\usepackage{url}
\usepackage[usenames]{color}
\usepackage[colorlinks=false]{hyperref}
\usepackage{cleveref}
\usepackage{tabularx}
\usepackage{xcolor}
\usepackage{mathabx}
\usepackage{enumitem}
\usepackage{fancyvrb}

\theoremstyle{plain}

\newtheorem{corollary}{Corollary}[section]
\newtheorem{lemma}{Lemma}[section]
\newtheorem{theorem}{Theorem}[section]

\newtheorem{problem}{Problem}[section]
\newtheorem{example}{Example}[section]

\theoremstyle{definition}

\crefname{conjecture}{Conjecture}{Conjectures}
\crefname{theorem}{Theorem}{Theorems}
\crefname{corollary}{Corollary}{Corollaries}
\crefname{lemma}{Lemma}{Lemmas}
\crefname{proposition}{Proposition}{Propositions}
\crefname{remark}{Remark}{Remarks}
\crefname{notation}{Notation}{Notations}
\crefname{problem}{Problem}{Problems}
\crefname{example}{Example}{Examples}
\crefname{section}{\S}{Sections}
\crefname{definition}{Definition}{Definitions}

\newcommand{\floor}[1]{\left\lfloor #1 \right\rfloor}

\title{Arithmetic Terms for Multinomial Coefficient Sums}
\author{Joseph M. Shunia \footnote{Wraithwatch, Austin, Texas, United States. E-mail address: {\tt jshunia@gmail.com}.}, Lorenzo Sauras-Altuzarra \footnote{Kurt Gödel Society, Vienna, Austria. E-mail address: {\tt lorenzo@logic.at}.}}
\date{May 2025}

\begin{document}

\maketitle

\begin{abstract} \noindent We construct arithmetic terms representing the partial sums of binomial coefficients, and we extend these results to obtain arithmetic terms representing the multisections of binomial coefficient sums. We also introduce an arithmetic term representing a certain type of multinomial coefficient sum and, as an application, we provide an arithmetic term representing the central trinomial coefficients. This solves one of the research problems of the celebrated book ``\textit{Concrete Mathematics}'', which remained open for nearly thirty years. \\[2 mm]
 \end{abstract}

\section{Introduction}

By abstracting the process of digit extraction (\cref{section:extraction}), present explicit combinatorial formulas, focusing on partial sums of binomial coefficients sums (\cref{section:partialsums}), multisections of binomial coefficient sums (\cref{section:multisectionsums}) and polynomial coefficients (\cref{section:polynomialcoefficients}). Our formulas for the partial sums of binomial coefficients appear to be the first of their kind, and it was assumed by some authors that such formulas cannot exist. The formula for the multisection sums of binomial coefficients remains fixed regardless of the choice of input parameters, contrasting with existing formulas of variable length. And the application of our polynomial coefficient formula resolves an open problem posed by Graham et al.\ \cite[Exercise 7.56]{GrahamEtAl} concerning the existence of a simple closed form expression for the central trinomial coefficients, which stood for nearly thirty years.

These formulas, with one exception (viz., the formula from \cref{theorem:partialsumsviaboardman}), are arithmetic terms, expressions whose study dates back to Robinson's foundational work in the 1950s on computability theory. Early research centered on broad theoretical questions, such as: \textit{what functions can be represented using only these elementary arithmetic operations?} (cf.\ Grzegorczyk \cite{Grzegorczyk}, Herman \cite{Herman}, and Robinson \cite{Robinson}). Grzegorczyk's work established that the class of the so-called Kalmar elementary functions is exactly $ \mathcal{E}^3 $ in the Grzegorczyk hierarchy, a classification that organizes primitive recursive functions by complexity (see Mazzanti \cite[Introduction]{Mazzanti} and Marchenkov \cite[Introduction]{Marchenkov}). Later, Mazzanti \cite{Mazzanti} proved that arithmetic terms are sufficiently expressive to generate the entire class of Kalmar elementary functions. Despite these advances, explicit constructions of arithmetic terms, particularly for combinatorial and number-theoretic applications, remained relatively underexplored until the recent works of Prunescu \cite{Prunescu} and Prunescu \& Sauras-Altuzarra \cite{PrunescuSaurasAltuzarra, PrunescuSaurasAltuzarra2, PrunescuSaurasAltuzarra3}.

Obtaining arithmetic terms for given Kalmar elementary functions is currently a challenging problem. While the aforementioned results by Mazzanti imply the existence of arithmetic terms for many prominent combinatorial and number-theoretic functions (see Mazzanti \cite[Introduction]{Mazzanti}), including the $n$-th prime number $p_n$ and the prime counting function $\pi(n)$, explicit constructions remain elusive. By developing and analyzing new arithmetic terms, we hope to provide new insights into combinatorics, number theory, and related fields.

\section{Digit extraction} \label{section:extraction}

The \textbf{truncated subtraction}, denoted by $ \dotdiv $, is the binary operation given by $ a \dotdiv b = \max ( a - b , 0 ) $ (see Vereshchagin \& Shen \cite[p.\ 141]{VereshchaginShen}). For simplicity, we may use $-$ in place of $\dotdiv$ when it is clear that $a-b \geq 0$.

Given a non-empty set $ V $ of variables, an \textbf{arithmetic term} is an element of the smallest set $ S $ of expressions such that:
\begin{enumerate}[noitemsep, topsep = 0pt]
    \item every non-negative integer is an element of $ S $,
    \item $ V \subseteq S $, and
    \item for every two elements $ a $ and $ b $ of $ S $, the expressions $ a + b $, $ a \dotdiv b $, $ a \cdot b $, $ \floor{a/b} $, $ a \bmod b $, and $ a^b $ are elements of $ S $
\end{enumerate}
(cf.\ Prunescu \cite{Prunescu} and Prunescu \& Sauras-Altuzarra \cite{PrunescuSaurasAltuzarra, PrunescuSaurasAltuzarra2, PrunescuSaurasAltuzarra3}).

\cref{example:marchenkov} is due to Marchenkov \cite[Section 2]{Marchenkov}.

\begin{example} \label{example:marchenkov} The values of the arithmetic terms $ m \cdot n $ and $$ \floor{\frac{2^{m + n + 4}}{\floor{\floor{2^{m + n + 4} / ( n + 1 )} / ( m + 1 )}}} \dotdiv ( m + n + 1 ) $$ coincide for every evaluation their arguments, $ m $ and $ n $. \end{example}

We follow the conventions $ 0^0 = 1 $ (see Mendelson \cite[Proposition 3.16]{Mendelson}) and $ \left\lfloor x / 0 \right\rfloor = 0 $ (see Mazzanti \cite[Section 2.1]{Mazzanti}), so, as the arithmetic term $ m \bmod n $ can be written as $ m \dotdiv ( n \cdot \floor{m / n} ) $, we also have that $ m \bmod 0 = m $ and $ m \bmod 1 = 0 $. In addition, in order to make the notation less cumbersome, we write the usual subtraction instead of the truncated subtraction in our formulas: they are interchangeable in all of them, with the exception of \cref{theorem:partialsumsviaboardman} (whose formula is not an arithmetic term).

\cref{lemma:wang}, which was contributed by Jinyuan Wang (pers.\ comm.), will be of utility in several occasions.

\begin{lemma} \label{lemma:wang}
Given two integers $ x \geq 2 $ and $ k \geq 1 $, and $ k $ non-negative integers $ a_0 $, $ \ldots $, $ a_{k - 1} $, each of them strictly smaller than $ x $, we have that
\begin{align*}
a_0 x^0 + \cdots + a_{k - 1} x^{k - 1} < x^k .
\end{align*}
\end{lemma}
\begin{proof} The numbers $ a_0 $, $ \ldots $, $ a_{k - 1} $ are integers, so $ a_i \leq x - 1 $ for every $ i \in \{ 0 , \ldots, k - 1 \} $ and thus
\begin{align*}
& a_0 x^0 + \cdots + a_{k - 1} x^{k - 1}
\leq ( x - 1 ) ( x^0 + \cdots + x^{k - 1} )
= ( x - 1 ) \frac{x^k-1}{x-1} 
= x^{k} - 1
< x^k .
\end{align*} \end{proof}

Given an integer $ k \geq 1 $, the $ k $-th digit of a given integer $ n \geq 0 $ in base $ b \geq 2 $, reading from right to left, is $ \floor{n / 10^{k - 1}} \bmod 10 $. For instance, $ \floor{7389056 / 10^{2 - 1}} \bmod 10 = 5 $ and $ \floor{738 / 2^{2 - 1}} \bmod 2 = 1 $ ($ 738 $ is written as $ 1011100010 $ in binary base). This process of digit extraction relies on \cref{theorem:extraction}, a property that we will utilize to express the coefficients of some important polynomials as arithmetic terms. Recall that, given a polynomial (or, with more generality, a formal power series) $ f ( x ) $ and a non-negative integer $ i $, the $ i $-th coefficient of $ f ( x ) $ is denoted by $ [ x^i ] f ( x ) $.

\begin{theorem} \label{theorem:extraction} If $ f ( x ) $ is a polynomial of non-negative integer coefficients, $ k \in \{ 0 , \ldots , \deg ( f ( x ) ) \} $ and $ c $ is an integer such that $ [ x^i ] f ( x ) < c $ for every $ i \in \{ 0 , \ldots , k \} $, then \begin{align*} [x^k]f(x) = \floor{\frac{f(c)}{c^k}} \bmod c . \end{align*} \end{theorem}

\begin{proof} Let $ r = \deg ( f ( x ) ) $ and, for every $ i \in \{ 0 , \ldots , r \} $, let $ a_i = [ x^i ] f ( x ) $.

If $ k = 0 $, then, as $ a_k < c $, we get $ \floor{f ( c ) / c^k} \bmod c = f ( c ) \bmod c = \sum_{i = 0}^r ( a_i c^i ) \bmod c = a_0 \bmod c = a_k $.

If $ c = 1 $, then $ a_k = 0 = \floor{f ( 1 ) / 1^k} \bmod 1 = \floor{f ( c ) / c^k} \bmod c $.

Now, suppose that $ k \geq 1 $ and $ c \geq 2 $. As we have, in addition, that $ a_i < c $ for every $ i \in \{ 0 , \ldots , k - 1 \} $, it follows, by applying \cref{lemma:wang}, that $ \sum_{i = 0}^{k - 1} a_i c^i < c^k $. Thus, $ 0 \leq \sum_{i = 0}^{k - 1} a_i c^{i - k} < 1 $. Moreover, it is clear that $ f ( c ) = \sum_{i = 0}^r a_i c^i $, so $ f ( c ) / c^k = \sum_{i = 0}^r a_i c^{i - k} $ and hence $ \floor{f ( c ) / c^k} = \sum_{i = k}^r a_i c^{i - k} $. Therefore, as $ a_k < c $, we get $ \floor{f ( c ) / c^k} \bmod c = \sum_{i = k}^r ( a_i c^{i - k} ) \bmod c = a_k \bmod c = a_k $. \end{proof}

\begin{example}
Consider the following polynomial:
\begin{align*}
    f(x) := 16x^6 + 24x^5 + 25x^4 + 20x^3 + 10x^2 + 4x + 1 .
\end{align*}
Applying \cref{theorem:extraction} to recover the coefficient $[x^2]f(x)$, one computes:
\begin{align*}
[x^2]f(x)
= \floor{\frac{f(100)}{100^2}} \bmod 100
= \floor{\frac{16242520100401}{10000}} \bmod 100
= 10 .
\end{align*}
One may repeat the process with $k \in \{ 0, 1, 3, 4, 5, 6 \}$ to recover the remaining coefficients.
\end{example}

An interesting observation underlying \cref{theorem:extraction} is that, if $ f ( x ) $ has more than one nonzero coefficient, then the choice $ c := f(1) $ (which is the sum of all the coefficients of $ f ( x ) $) serves as a suitable base to recover all the coefficients of $ f(x) $.

It is worth mentioning that the application of \cref{theorem:extraction} is closely related to an algorithmic technique known as \textbf{Kronecker substitution} (see von zur Gathen and Gerhard \cite[Section 8.4]{vonzurGathenandGerhard}). This technique encodes a polynomial $ f(x) $ of non-negative integer coefficients as an integer, by selecting a suitable base $ d \geq 2 $ such that the value $ f(d) $ encodes the coefficients of $ f(x) $ in its base-$ d $ digit expansion. Kronecker substitution has been important in the development of fast polynomial multiplication algorithms (see Albrecht et al.\ \cite{AlbrechtEtAl}, Bos et al.\ \cite{BosEtAl}, Greuet et al.\ \cite{GreuetEtAl}, Harvey \cite{Harvey}, and Harvey \& van der Hoeven \cite{HarveyandvanderHoeven}).

To demonstrate the utility of \cref{theorem:extraction} in combinatorial and number-theoretic contexts, we present a new proof of a binomial coefficient identity that was first established by Robinson \cite{Robinson}.

\begin{corollary}
If $ k $ and $ n $ are integers such that $ n \geq 1 $ and $ 0 \leq k \leq n $, then
\begin{align*}
\binom{n}{k} = \floor{\frac{(2^n+1)^n}{2^{n k}}} \bmod 2^n .
\end{align*}
\end{corollary}
\begin{proof} By applying the binomial theorem, we have that $ f(x) := ( x + 1 )^n = \sum_{j=0}^n \binom{n}{j} x^j $, whence it follows that $ 0 \leq [ x^i ] f ( x ) = \binom{n}{i} < \sum_{j=0}^n \binom{n}{j} = ( 1 + 1 )^n = 2^n $ for every $ i \in \{ 0 , \ldots , k \} $. Therefore, by applying \cref{theorem:extraction}, we get $ \binom{n}{k} = \floor{f(2^n) / (2^n)^k} \bmod 2^n $, whence the statement immediately follows. \end{proof}

\section{Partial sums of binomial coefficients} \label{section:partialsums}

We say that the \textbf{consecutive term ratio} of an expression $ e ( n ) $ is the expression $ e ( n + 1 ) / e ( n ) $. In addition, a \textbf{hypergeometric term} with respect to a field $ K $ is defined as a univariate expression whose consecutive term ratio is a rational function on $ K $ (i.e.\ an expression $ e ( n ) $ such that $ e ( n + 1 ) / e ( n ) \in K ( n ) $). And a \textbf{hypergeometric closed form} with respect to $ K $, as a linear combination of hypergeometric terms. See Petkovšek et al.\ \cite[Definition 8.1.1]{PetkovsekEtAl} and Sauras-Altuzarra \cite[Definition 1.4.13]{SaurasAltuzarra}.

Boardman \cite{Boardman} asserted ``\textit{it is well-known that there is no closed form (that is, direct formula) for the partial sum of binomial coefficients}'' and referred to Petkovšek et al.\ \cite{PetkovsekEtAl}. However, Petkovšek et al.\ \cite[p.\ 88 and 102]{PetkovsekEtAl} only mention that it is not representable as a hypergeometric term, it is not clear from their book whether it has a hypergeometric closed form or not. Anyway, we now show how to represent it as an arithmetic term, in two different ways.

\begin{theorem} \label{theorem:partialsums}
If $ j $ and $ n $ are integers such that $ 0 \leq j \leq n - 2 $, then $$ \sum_{k=0}^{j} \binom{n}{k} = \floor{\frac{(2^n+1)^n}{2^{n(n-j)}}} \bmod (2^n-1) . $$
\end{theorem}

\begin{proof} The first step in formula (i) is to perform floored division by $2^{n(n-j)}$ on each side of the equation $(2^n+1)^n = \sum_{k=0}^n \binom{n}{k} 2^{nk}$, which yields \begin{equation} \label{equation:partialsums} \frac{( 2^n + 1 )^n}{2^{n ( n - j )}} = \sum_{k=0}^n \binom{n}{k} 2^{nk-n(n-j)} . \end{equation}

We have that $ 0 \leq j < n $, so $ 2^n \geq 2 $, $ n - j \geq 1 $, and $ \binom{n}{k} < \sum_{i = 0}^n \binom{n}{i} = ( 1 + 1 )^n = 2^n $ for every $ k \in \{ 0 , \ldots , n - j - 1 \} $.

Then, by applying \cref{lemma:wang}, we have that $ \sum_{k = 0}^{n - j - 1} \binom{n}{k} ( 2^n )^k < ( 2^n )^{n - j} $. That is to say, we get $ \sum_{k = 0}^{n - j - 1} \binom{n}{k} 2^{nk-n(n-j)} < 1 $, whence it follows, by \cref{equation:partialsums} and the fact that $ [ n k - n ( n - j ) \geq 0 \ \Leftrightarrow \ k \geq n - j ] $, that $$ \floor{\frac{(2^n+1)^n}{2^{n(n-j)}}} = \sum_{k=n-j}^n \binom{n}{k} 2^{nk-n(n-j)} . $$

The last equation can be rewritten as $$ \floor{\frac{(2^n+1)^n}{2^{n(n-j)}}} = \sum_{k=0}^j \binom{n-k}{k} 2^{n(j-k)} $$ and, because of the symmetry for binomial coefficients in row $n$ (i.e.\ $\binom{n}{k} = \binom{n}{n-k}$), it reduces to $$ \floor{\frac{(2^n+1)^n}{2^{n(n-j)}}} = \sum_{k=0}^j \binom{n}{k} 2^{n(j-k)} . $$ Therefore, by reducing modulo $ 2^n - 1 $ (which replaces all instances of $2^n$ with one, by the polynomial remainder theorem and the fact that $ \sum_{k=0}^j \binom{n}{k} < \sum_{k=0}^{n - 1} \binom{n}{k} = \sum_{k=0}^n \binom{n}{k} - \binom{n}{n} = ( 1 + 1 )^n - 1 = 2^n - 1 $), the statement follows. \end{proof}

We now provide an alternative formula for the partial sums of binomial coefficients, using results from Boardman \cite{Boardman}.

\begin{theorem} \label{theorem:partialsumsviaboardman}
If $ j $ and $ n $ are integers such that $ 1 \leq j \leq n $, then
\begin{align*}
\sum_{k=0}^{j} \binom{n}{k} = 1 + \left( \floor{\frac{1 - ( 2^n + 1 )^n}{( 2^n - 1 ) 2^{n j}}} \bmod 2^n \right) .
\end{align*}
\end{theorem}

\begin{proof} Let $ x $ be a complex number $ x $ such that $ | x | < 1 $.

Boardman \cite{Boardman} proved that $$ \sum_{u = 1}^\infty \sum_{v = 1}^u \binom{n}{v} x^u = \frac{( x + 1 )^n - 1}{1 - x} . $$

As $ \sum_{i = 0}^n \binom{n}{i} = ( 1 + 1 )^n = 2^n $ and $ \binom{n}{k} $ for every integer $ k > n $, we have that $$ \sum_{u = n + 1}^\infty \sum_{v = 1}^u \binom{n}{v} x^u = \sum_{u = n + 1}^\infty ( 2^n - 1 ) x^u = ( 2^n - 1 ) x^{n + 1} \sum_{u = 0}^\infty x^u = \frac{( 2^n - 1 ) x^{n + 1}}{1 - x} . $$

Hence, we have that $$ f ( x ) := \sum_{u = 1}^n \sum_{v = 1}^u \binom{n}{v} x^u = \frac{- ( 2^n - 1 ) x^{n + 1} + ( x + 1 )^n - 1}{1 - x} = \frac{( 2^n - 1 ) x^{n + 1} - ( x + 1 )^n + 1}{x - 1} . $$

In addition, $ \sum_{k = 1}^i \binom{n}{k} < \sum_{k = 0}^n \binom{n}{k} = ( 1 + 1 )^n = 2^n $ for every $ i \in \{ 1 , \ldots , j \} $.

Therefore, by applying \cref{theorem:extraction}, we have that $$ \sum_{k=1}^{j} \binom{n}{k} = [ x^j ] f ( x ) = \floor{\frac{( 2^n - 1 ) ( 2^n )^{n + 1} - ( 2^n + 1 )^n + 1}{( 2^n - 1 ) ( 2^n )^j}} \bmod 2^n = $$ $$ \floor{2^{n ( n + 1 - j )} + \frac{1 - ( 2^n + 1 )^n}{( 2^n - 1 ) 2^{n j}}} \bmod 2^n = \floor{\frac{1 - ( 2^n + 1 )^n}{( 2^n - 1 ) 2^{n j}}} \bmod 2^n , $$ whence the statement immediately follows. \end{proof}

Note that the right-hand side of the identity from \cref{theorem:partialsumsviaboardman} is not an arithmetic term, because the expression $ 1 - ( 2^n + 1 )^n $ cannot be replaced with $ 1 \dotdiv ( 2^n + 1 )^n $.

\section{Multisections of binomial coefficient sums} \label{section:multisectionsums}

An \textbf{arithmetic progression} is a sequence of the form $ ( k s + j )_{k = 0}^\infty $, where $ s $ and $ j $ are integers such that $ s \neq 0 $ (cf.\ Weisstein \cite{Weisstein}). And the \textbf{multisection} of a formal power series $ \sum_{k = 0}^\infty f ( k ) x^k $, with respect to an arithmetic progression $ a $ of non-negative terms, is the formal power series $ \sum_{k = 0}^\infty f ( a ( k ) ) x^{a ( k )} $ (cf.\ Comtet \cite[p.\ 84]{Comtet}). In this section we focus on the case in which $ f ( k ) = \binom{n}{k} $, where $ n $ is some positive integer, and $ x $ is one; that is to say, on the binomial coefficient sums. A classic formula is \begin{equation} \label{equation:traditionalmultisections} \sum_{v=0}^\infty \binom{n}{vs+j} = \frac{1}{s} \sum_{k=0}^{s-1} \left(2 \cos\left(\frac{\pi k}{s}\right) \right)^n \cos\left( \frac{\pi (n-2j) k}{s} \right) \end{equation} (recall that $ \binom{n}{i} = 0 $ for every integer $ i > n $), which holds for every three integers $ j $, $ n $, and $ s $ such that $ n \geq 1 $ and $ 0 \leq j < s $, and whose length depends on $ s $ (cf.\ Weisstein \cite{Weisstein2}). In contrast, we present an arithmetic term for this case, providing a significant simplification of the computation.

\begin{theorem} \label{theorem:multisections} Given integers $ j \geq 0 $, $ n \geq 1 $, and $ s \geq 2 $ such that $ j < s $ and $ ( s - 1 ) n \geq 2 $, we have that \begin{align*} \sum_{v=0}^{\floor{(n-j)/s}} \binom{n}{vs+j} &= \floor{\frac{(2^n+1)^n \bmod (2^{ns}-1)}{2^{nj}}} \bmod 2^n . \end{align*} \end{theorem}

\begin{proof} By the binomial theorem, the number $ ( 2^n + 1 )^n \bmod ( 2^{n s} - 1 ) $ is equal to $$ \left( \sum_{k = 0}^{n} \binom{n}{k} 2^{k n} \right) \bmod ( 2^{n s} - 1 ) . $$ Observe that, given two non-negative integers $ u $ and $ v $, the condition $ v s + u \leq n $ can be rewritten as $ v \leq ( n - u ) / s $ or, equivalently, as $ v \leq \floor{( n - u ) / s} $. Hence, we can rewrite the previous sum as $$ \left( \sum_{u = 0}^{s - 1} \sum_{v = 0}^{\floor{( n - u ) / s}} \binom{n}{v s + u} ( 2^{v s + u} )^n \right) \bmod ( 2^{n s} - 1 ) . $$ Thus, by reducing modulo $ 2^{n s} - 1 $ (which replaces all instances of $ 2^{n s} $ with one, by the fact that $$ \sum_{u = 0}^{s - 1} \sum_{v = 0}^{\floor{( n - u ) / s}} \binom{n}{v s + u} 2^{u n} = \sum_{r = 0}^n \binom{n}{r} 2^{( r \bmod s ) n} = 2^{n s} - \left( 2^{n s} - \sum_{r = 0}^n \binom{n}{r} 2^{( r \bmod s ) n} \right) = $$ $$ 2^{n s} - \left( 2^{( s - 1 ) n} \sum_{r = 0}^n \binom{n}{r} - \sum_{r = 0}^n \binom{n}{r} 2^{( r \bmod s ) n} \right) = 2^{n s} - \sum_{r = 0}^n \binom{n}{r} ( 2^{( s - 1 ) n} - 2^{( r \bmod s ) n} ) \leq $$ $$ 2^{n s} - \binom{n}{0} ( 2^{( s - 1 ) n} - 2^{( 0 \bmod s ) n} ) = 2^{n s} - ( 2^{( s - 1 ) n} - 1 ) < 2^{n s} - 1 $$ and the polynomial remainder theorem), the previous sum becomes $$ \sum_{u = 0}^{s - 1} \sum_{v = 0}^{\floor{( n - u ) / s}} \binom{n}{v s + u} 2^{u n} $$ and therefore \begin{equation} \label{equation:multisections} \frac{( 2^n + 1 )^n \bmod ( 2^{n s} - 1 )}{2^{n j}} = \sum_{u = 0}^{s - 1} \sum_{v = 0}^{\floor{(n-u)/s}} \binom{n}{v s + u} 2^{n ( u - j )} . \end{equation}

In addition, we know that $ n \geq 1 $ and $ s \geq 2 $, so, for every $ u \in \{ 0 , \ldots , j - 1 \} $, we get $$ \sum_{v = 0}^{\floor{(n-u)/s}} \binom{n}{v s + u} < \sum_{v = 0}^n \binom{n}{v} = ( 1 + 1 )^n = 2^n . $$ Then, by applying \cref{lemma:wang}, we have that $$ \sum_{u = 0}^{j - 1} \sum_{v = 0}^{\floor{(n-u)/s}} \binom{n}{v s + u} ( 2^n )^u < ( 2^n )^j . $$ And this inequality is equivalent to $$ \sum_{u = 0}^{j - 1} \sum_{v = 0}^{\floor{(n-u)/s}} \binom{n}{v s + u} 2^{n ( u - j )} < 1 , $$ whence it follows, by applying that $ j < s $ and \cref{equation:multisections}, that $$ \floor{\frac{( 2^n + 1 )^n \bmod ( 2^{n s} - 1 )}{2^{n j}}} = \sum_{u = j}^{s - 1} \sum_{v = 0}^{\floor{(n-u)/s}} \binom{n}{v s + u} 2^{n ( u - j )} . $$ The statement is now straightforward. \end{proof}

\section{Polynomial coefficients} \label{section:polynomialcoefficients}

A number of the form $ \frac{( k_1 + \cdots + k_r ) !}{k_1 ! \cdot \ldots \cdot k_r !} $, where $ r $ is a positive integer and $ k_1 $, $ \ldots $, $ k_r $ are non-negative integers, is known as a \textbf{multinomial coefficient} and denoted by $ \binom{k_1 + \cdots + k_r}{k_1 , \ldots, k_r} $ (see Brualdi \cite[Section 5.4]{Brualdi}). In fact, the so-called \textbf{multinomial theorem} asserts that this number is the coefficient of the monomial $ x_1^{k_1} \cdot \ldots \cdot x_r^{k_r} $ in the polynomial $ ( x_1 + \cdots + x_r )^{k_1 + \cdots + k_r} $ (cf.\ Brualdi \cite[Theorem 5.4.1]{Brualdi}).

Given three integers $ k \geq 0 $, $ n \geq 0 $ and $ r \geq 1 $ such that $ k \leq n ( r - 1 ) $, the number $ [ x^k ] ( x^0 + \cdots + x^{r - 1} )^n $ is called \textbf{polynomial coefficient} and denoted by $ \binom{n , r}{k} $ (cf.\ Comtet \cite[p.\ 77]{Comtet}). One way of computing this number consists of applying the identity \begin{equation*} \binom{n , r}{k} = \frac{1}{2 \pi} \int_0^{2 \pi} \left( \frac{\exp ( r \iota x ) - 1}{\exp ( \iota x ) - 1} \right)^n \exp ( - k \iota x ) dx \end{equation*} (cf.\ André \cite[Corollaire 70]{Andre} (in French)), where $ \iota $ denotes the imaginary unit (a notation used, for instance, by Schneider \cite{Schneider}) or, equivalently, \begin{equation*} \binom{n , r}{k} = \frac{2}{\pi} \int_0^{\pi / 2} \left( \frac{\sin ( r x )}{\sin ( x )} \right)^n \cos ( ( n ( r - 1 ) - 2 k ) x ) dx \end{equation*} (cf.\ Li \cite[p.\ 3]{Li}).

Rudolph-Lilith \& Muller \cite[Lemma 2]{RudolphLilithMuller} obtained the identity $$ \binom{n , 2 k + 1}{k n} = \frac{1}{2 k n + 1} \left( ( 2 k + 1 )^n + \sum_{j = 1}^{2 k n} \left( \frac{\sin ( ( 2 k + 1 ) j \pi / ( 2 k n + 1 ) )}{\sin ( j \pi / ( 2 k n + 1 ) )} \right)^n \right) , $$ where $ k $ and $ n $ are non-negative integers. Unfortunately, as in the case of \cref{equation:traditionalmultisections}, this formula is of variable length. However, by means of \cref{theorem:extraction}, we can now express the polynomial coefficients as arithmetic terms.

\begin{theorem} \label{theorem:polynomialcoefficients} If $ k $, $ n $ and $ r $ are integers such that $ n \geq 1 $, $ r \geq 2 $, and $ 0 \leq k \leq n ( r - 1 ) $, then
\begin{align*}
    \binom{n , r}{k} = \floor{\left(\frac{r^{rn} - 1}{r^{n+k} - r^k}\right)^n} \bmod r^n .
\end{align*}
\end{theorem}

\begin{proof} Let $ f ( x ) = ( x^0 + \cdots + x^{r - 1} )^n $, which is a polynomial of positive integer coefficients such that $ \deg ( f ( x ) ) = n ( r - 1 ) $.

We know that $ n \geq 1 $ and $ r \geq 2 $, so $ n ( r - 1 ) \geq 1 $ and then, for every $ i \in \{ 0 , \ldots , k \} $, we have that $$ [ x^i ] f ( x ) = \binom{n , r}{i} < \sum_{u = 0}^{n ( r - 1 )} \binom{n , r}{u} = \sum_{u = 0}^{n ( r - 1 )} [ x^u ] f ( x ) = f ( 1 ) = r^n . $$

Therefore, by applying \cref{theorem:extraction} and the well-known fact that $ f ( x ) = ( ( x^r - 1 ) / ( x - 1 ) )^n $, we get $$ \binom{n , r}{k} = \floor{\frac{f ( f ( 1 ) )}{f ( 1 )^k}} \bmod f ( 1 ) = \floor{\left( \frac{( r^n )^r - 1}{r^n - 1} \right)^n \frac{1}{( r^n )^k}} \bmod r^n = \floor{\left(\frac{r^{rn} - 1}{r^{n+k} - r^k}\right)^n} \bmod r^n . $$ \end{proof}

We now obtain a formula for the partial sums of polynomial coefficients, by taking a similar approach as in \cref{section:partialsums}.

\begin{theorem} \label{theorem:polynomialcoefficientssums} If $ j $, $ n $ and $ r $ are integers such that $ n \geq 1 $, $ r \geq 2 $, and $ 0 \leq j \leq n ( r - 1 ) - 2 $, then
\begin{align*}
\sum_{k=0}^{j} \binom{n , r}{k}
= \left( \left( \frac{r^{rn}-1}{r^n-1} \right)^n \bmod r^{n(j+1)} \right) \bmod (r^n-1) .
\end{align*}
\end{theorem}

\begin{proof} Let $ f ( x ) $ be the polynomial $ ( x^0 + \cdots + x^{r - 1} )^n $, which is equal to $ ( ( x^r - 1 ) / ( x - 1 ) )^n $.

It is clear that $$ \sum_{k=0}^{n ( r - 1 )} \binom{n , r}{k} (r^n)^k = \sum_{k=0}^{n ( r - 1 )} [ x^k ] f ( x ) (r^n)^k = f ( r^n ) = \left( \frac{r^{nr}-1}{r^n-1} \right)^n . $$

In addition, for every $ k \in \{ 0 , \ldots , j \} $, we have that $$ \binom{n , r}{k} < \sum_{u = 0}^{n ( r - 1 )} \binom{n , r}{u} = \sum_{u = 0}^{n ( r - 1 )} [ x^u ] f ( x ) = f ( 1 ) = r^n . $$ Thus, by applying \cref{lemma:wang}, we have that $$ \sum_{k=0}^j \binom{n , r}{k} r^{n k} < r^{n ( j + 1 )} $$ and, consequently, that $$ \sum_{k=0}^j \binom{n , r}{k} r^{nk} = \left( \frac{r^{nr}-1}{r^n-1} \right)^n \bmod r^{n ( j + 1 )} . $$

Thus, by reducing modulo $ r^n - 1 $ (which replaces all instances of $ r^n $ with one, by the fact that $$ \sum_{k=0}^{j} \binom{n , r}{k} < \sum_{k = 0}^{n ( r - 1 ) - 1} \binom{n , r}{k} = \sum_{k = 0}^{n ( r - 1 )} \binom{n , r}{k} - \binom{n , r}{n ( r - 1 )} = $$ $$ \sum_{k = 0}^{n ( r - 1 )} [ x^k ] f ( x ) - [ x^{n ( r - 1 )} ] f ( x ) = f ( 1 ) - 1 = r^n - 1 $$ and the polynomial remainder theorem), the statement follows. \end{proof}

The expressions $ \binom{n , 3}{n} $, often called \textbf{central trinomial coefficients}, do not have a hypergeometric closed form (see Petkovšek et al.\ \cite[Theorem 8.8.1]{PetkovsekEtAl}). Based on this result, Graham et al.\ \cite[Exercise 7.56]{GrahamEtAl} posed the following research problem.

\begin{problem} \label{problem:grahametal} Prove that there is no simple closed form for the coefficient of $x^n$ in $(x^2+x+1)^n$, as a function of $n$, in some large class of simple closed forms. \end{problem}

Here, a ``simple closed form'' is defined as an expression using only addition, subtraction, multiplication, division, and exponentiation, in explicit ways (see Graham et al.\ \cite[p.\ 7]{GrahamEtAl}). This is essentially the definition of an arithmetic term. And, by applying our polynomial coefficient formula from \cref{theorem:polynomialcoefficients}, we see that
\begin{align*}
    \binom{n , 3}{n} = \floor{\left(\frac{3^{3n} - 1}{3^{2n} - 3^n}\right)^n} \bmod 3^n
    = \floor{\left(\frac{27^n - 1}{9^n - 3^n}\right)^n} \bmod 3^n
\end{align*}
for every integer $ n \geq 1 $, which provides a negative answer to \cref{problem:grahametal}.

\end{document}